**А.Я. Белов**, г. Рамат-Ган, Израиль
**Р. Явич**, г. Нетания, Израиль

# Проблемы одарённости и стадийность математического обучения (к работе И.С. Рубанова «Как обучать методу математической индукции»).

## Проблема неправильной личностной ориентации

Проблема одаренности, помимо традиционного, имеет и другой аспект. Так ли «бездарны» остальные дети? Не является ли в ряде случаев т.н. «бездарность» следствием неправильной ориентации личности? Не является ли эта ориентация следствием некоторых аспектов традиционной системы образования? И если да, то, как достичь правильной ориентации?

Представляется необходимым сравнить не только «интеллектуальные качества» «одаренных» и «бесталанных» учащихся, но и поведенческие особенности применительно к предмету. (Говоря о т.н. «супердостижениях» также следует отметить, что они зависят как от интеллектуальных качеств биологического тела, так и от поведенческих особенностей.)

В процессе преподавания возникает следующее наблюдение. «Слабые» учащиеся относятся к предмету, как к своду неких правил и инструкций. И главное — не нарушить правило, не совершить ошибку. Те учащиеся (школьники или студенты), которые смогли преодолеть такое отношение к предмету, учатся значительно эффективнее и называются *одаренными*. Спору нет, у них уже есть заслуга. Когда говорят, что «студенты не понимают понятие скорости» — это и правда, и неправда. Студенты действительно затрудняются дать формальное определение производной, но интуитивное представление о скорости есть у любого нормального человека.

Точно так же, каждому нормальному человеку ясно, что с практической точки зрения яблоко бесконечно больше атома и бесконечно меньше Земли. В то же время абстрактные рассуждения с величинами разной малости усваиваются с большим трудом. Причина не в том, что эти рассуждения сложные, а в некотором культурном барьере.

Опыт говорит, что главное отличие настоящего отличника от рядового студента заключается, прежде всего, не столько в интеллекте, сколько в том, что у отличника здравый смысл направлен на сам предмет. Содержание занятий в подавляющем большинстве случаев особых требований к интеллекту не предъявляет (традиционная методика старается максимально уйти от трудных доказательств). Собственно, творческая одаренность проявляется на более высоком уровне.

## Стадийность образования и пути решения методической проблемы

**Опыт кружкового преподавания.** Есть мнение: студенты вуза «слабые» и их нельзя учить, так как учат 12-летних школьников на кружках. Методы кружкового преподавания годятся только для сильных. Это значит, что студентов надо учить традиционно, обучая стандартным приемам решения задач из стандартного материала. От этого они, впрочем, не становятся сильнее.

Кружки интересны, прежде всего, тем, что отбор материала более гибок чем в стандартной школьной или вузовской программе. Он управлялся интересами учащихся и преподавателей, он более *естественный*. С более разнообразным материалом происходит экспериментирование. Мы убеждены, что накопленный опыт может оказаться важным и для преподавания «слабым» учащимся.

**Стадийность математического образования.** При первом знакомстве математика — это, прежде всего, наука о решении занимательных задач и головоломок (и немного — о некоторых интересных природных закономерностях). Первоначальное обращение к математике должно быть именно таким узким и на первых порах следует избегать доказательства вещей, которые «и так очевидны». Вначале нужно увидеть, как строгие рассуждения позволяют получить нетривиальные результаты, и только потом осознать, что же такое *строгое рассуждение*. Но это — работа следующего этапа.

Изучение математики вообще и любой темы, в частности, состоит из нескольких этапов, которые нельзя смешивать. Первый этап заключается в формировании интуиции («мыслеобразов»), затем следует собственно *осознание* и затем — *формирование умения управлять своим мышлением*.

Как говорил Р.Бернс «Твой стиль, суховатый и сдержанно краткий, без умолку хвалят друзья. Уздечка нужна, чтобы править лошадкой, но где же лошадка твоя?». Традиционное преподавание начинается именно с уздечки.

Преждевременное наведение строгости крайне вредно. Оно вызывает т.н. «эффект сороконожки»[1].

**Обсуждение работы И.С. Рубанова.**

Формально посвященная методу математической индукции, эта работа на самом деле имеет гораздо более широкий методический и методологический интерес, и жанр комментариев к ней и размышлений вокруг нам представляется уместным.

Традиционное преподавание этого метода состоит в его словесном описании, формализации, манипуляции закрепляются упражнениями на доказательство тождеств.

Однако необходимо осуществить предварительное формирование *интуитивного представления* следуя *принципу наглядности*. Необходимо создать набор мыслеобразов, которые необходимы для размышления. Такими полезными мыслеобразами являются образ застежки-молнии, волны доказательств, ходьбы по лестнице, наведения мостов к более сложным задачам.

---

[1] Согласно легенде, сороконожка разучилась ходить, когда ее убедили сознательно двигать каждым членом. Отсюда возник соответствующий психологический термин.

Оформление метода индукции через «базу» и «переход» следует осуществить позднее. При этом все начинается с интуитивного образа: идя по цепочке, мы сможем дойти до любого утверждения. Далее обсуждается как сократить бесконечную цепочку утверждений – обсуждается сокращенная запись или *кодирование*. Происходит *осознание* бесконечной цепочки утверждений в виде единой теоремы.

Возникает понятие *переменной*. (Понятие *переменной* является чрезвычайно важным и нетривиальным. Упомянем методические проблемы связанные с преподаванием производной композиции функций. Возникающие проблемы связаны с понятием *предиката*. Все эти понятия обладают самоценностью помимо ценности самого метода индукции.)

**Ключевые задачи.**

Общая схема изучения метода математической индукции состоит в том, что сценарий проигрывается на нескольких ключевых задачах отражающих своего рода «абстрактные ядра» рассуждений. Вначале разбирается, например,

**Задача.** *Число составлено из 27 единиц, идущих подряд. Доказать, что оно делится на 27.*

**Идея доказательства** состоит в разбиении на 3 блока по 9 единиц идущих подряд. Каждый из них делится на 9, а произведение такого бока на число 1000000001000000001 делится на 27.

Далее разбирается следующая

**Задача.** *Число составлено из 81 единицы, идущих подряд. Доказать, что оно делится на 81.*

Затем идет работа с аналогичной задачей о числе, составленном из 243 единиц и т. д. В конце концов, учащимся надоедает решать задачи-близнецы и они захотят рассуждать «в общем виде». Но как сформулировать «общий вид» утверждения? С этого начинается *осознание*.

Как говорил Д.Пойя, метод — это идея, примененная дважды. Разбирается также еще несколько серий задач подобного рода. Например, такая

**Задача.** *4 прямые разбивают плоскость на области. Доказать, что можно покрасить ее в черный и белый цвет так, чтобы соседние области были раскрашены в разные цвета.*

**Идея решения** состоит в том, что разбиение плоскости тремя прямыми порождает раскраску для разбиения плоскости четырьмя прямыми. Далее обсуждается случай большего числа прямых и работа производится по той же схеме, что и для предыдущей серии. Точно так же у учащиеся формируется потребность *рассуждать в общем виде*.

Итак, **первый этап** в преподавании индукции заключается в разборе нескольких ключевых задач. При этом формируется интуитивное ощущение, что задача для чего-то большего часто сводится к задаче для чего-то меньшего.

**Формулировка** метода математической индукции получается только после того, как соответствующие рассуждения будут проведены в общем виде,

и на первых порах метод математической индукции будет назван *«планом решения задачи»*.

**Преподавание и ключевые задачи.**

Дадим для начала примерное описание того, что же такое «ключевая задача» применительно к методу математической индукции.

1. Это задача *оптимальной сложности*. Если задача слишком простая, то она выглядит недостаточно убедительно и не раскрываются рассуждения, если слишком сложная, то трудности самого предмета начинают отвлекать.
2. Основное содержание решения составляет применение индукции, нет сильно отвлекающих от него посторонних трудностей.
3. Хорошо иллюстрирует метод. В данном случае легко разворачиваются в цепочки *неочевидных* утверждений (иначе неубедительно — см. п. 1). По нескольким первым утверждениям естественно разворачивается вся цепочка.
4. Задача интересна ученикам (ср. с п.1).

Техника ключевых задач важна отнюдь не только для преподавания метода математической индукции, но и для других тем. В связи с этим следует дать менее специализированное описание этого понятия.

Ключевая задача — это задача, в процессе решения которой достаточно хорошо усваивается изучаемая идея или метод. Она характеризуется следующими признаками.

1. Оптимальная сложность.
2. Чистота при использовании и возникновения метода, отсутствие отвлекающих обстоятельств.
3. Интересность и содержательность, естественность и мотивированность. Создается своего рода искусственная ситуация рождения метода.

**Первый этап** состоит в интуитивном ощущении сути. Типичный сценарий работы — повторение на новых сериях ключевых задач. Для демонстрации того, что это действительно *метод*, необходимо несколько таких серий. Способ, которым осуществляется разбор, в частности степень подробности, зависит больше от ученика, чем от специфики метода.

Процесс обучения любому предмету, в том числе математике, начинается *с создания интуиции*. Это означает *создание поля мыслеобразов*, которыми обрабатывается тот или иной круг задач. Иногда говорят о *словоформах*.

Создание интуиции также происходит в два этапа. **Первый этап** состоит «в вызывательной процедуре» по выражению оккультистов, т.е. в привлечении мыслеобразов из здравого смысла, наглядных представлений. Механизм *вызова* или *привлечения* (мыслеформ, словоформ) особенно важен в преподавании математики нематематикам и, вообще, для взрослого человека, когда он хочет изучить тот или иной предмет. Эти «бытовые мыслеформы» служат исходным материалом при формировании (можно сказать — преподавании) мыслеобразов, относящихся к самому предмету.

Р. Фейнман в своих знаменитых лекциях по физике активно занимается вызовом мыслеобразов. Например, разбирая понятие *мгновенной скорости* он приводит диалог полицейского и дамы ( – Вы ехали со скоростью 60 миль в час. – Простите, Сэр, но я за последний час 60 миль не проехала. – Но Вы ехали со скоростью 30 метров в секунду. – Но штрафуют за 60 миль в час, а не за 30 метров в секунду).

В привлечении мыслеобразов и состоит так называемый «*принцип наглядного обучения*».

**Второй** этап состоит в формировании профессиональных мыслеформ.

**Следующий этап** состоит в переходе к осознанию, к управлению мыслеобразами. От интуитивного образа идем к абстрактному ядру, к схеме решения, в конце происходит формализация.

**Индукция revisiting, повторение.**

Повторение и закрепление материала должно основываться на его показе с другой стороны, в непривычном ракурсе.

Связь между методом математической индукции и догадкой по аналогии слабее, чем кажется изначально. Математическая индукция служит средством проверки предварительно угаданных закономерностей. Такого рода деятельность хороша при повторении, помогает осознать метод, разгрузить задачи первого этапа изучения.

**Доказательство тождеств, свойств делимости, неравенств, индукция и приращение.**

Хорошо для повторения, но отнюдь не для первого знакомства с методом. Такого рода доказательства, как отмечал И.С. Рубанов, предполагают весьма нетривиальное для начинающего и в то же время недостаточно осознаваемое продвинутым человеком, осознание формулы как цепочки утверждений. Это опять-таки связано с понятием предиката и т.д. В качестве разгрузки полезно нахождение формул общего члена, продолжения последовательностей. Упомянем также методические проблемы, связанные с очевидностью первых равенств, а также с понятием *подстановки*.

**Вариации на тему индукции.** Нам кажется, что в подходе И.С. Рубанова есть некая непоследовательность. Разветвленные типы индукции, по всей видимости, стоит давать одновременно со стандартным ее типом на этапе формирования мыслеобраза. Тем самым мы показываем несколько схем и способов организации мысли и разгружаем дальнейшие трудности, связанные с осознанием, демонстрируя содержательное разнообразие схем. (Хорошей иллюстрацией служит неравенство Коши, несмотря на избитость этого сюжета.) То же относится и к двойной индукции.

При преподавании принципа вполне упорядоченности продвинутым участникам, полезно упомянуть, что этот принцип есть обобщение метода математической индукции.

**Некоторые технологические наблюдения.** Контрпримеры позволяют обойти подводные камни, неявные подсказки помогают и в отвлечения внимания (от ненужного пути решения) а главное — создают ощущение самостоятельности. Очень важен принцип разгрузки и разделения педагогических целей.

**В заключение** отметим, что преподавание процесс творческий и что создание мыслеобразов неоднозначно, не может делаться по шаблону и зависит от мастерства преподавателя.

**Литература.**